\documentclass[letterpaper, 10 pt, conference]{ieeeconf}  

\IEEEoverridecommandlockouts                              

\usepackage{graphicx,graphics,color,epsfig}
\usepackage{cite}
\usepackage{hyperref}
\usepackage{amssymb}  
\usepackage{tikz}
\usepackage{amsmath}
\usepackage{textcomp}
\usepackage{fixltx2e}
\usepackage{stfloats}
\usepackage{url}
\usepackage{framed}
\usepackage{pgfplots}
\usepackage[utf8]{inputenc}

\newenvironment{customlegend}[1][]{%
    \begingroup
    \csname pgfplots@init@cleared@structures\endcsname
    \pgfplotsset{#1}%
}{%
    \csname pgfplots@createlegend\endcsname
    \endgroup
}%

\def\addlegendimage{\csname pgfplots@addlegendimage\endcsname}

\pgfkeys{/pgfplots/number in legend/.style={%
        /pgfplots/legend image code/.code={%
            \node at (0.125,-0.0225){#1}; 
        },%
    },
}
\pgfplotsset{
every legend to name picture/.style={west}
}

\title{\LARGE \bf 
  A discrete event traffic model explaining the traffic phases of the train dynamics in a metro line system with a junction}

\author{Florian Schanzenbächer$^{1*}$, Nadir Farhi$^{2}$, Zoi Christoforou$^{3}$, Fabien Leurent$^{3}$, Gérard Gabriel$^{1}$
\thanks{$^{1}$ RATP, Paris, France.}%
\thanks{$^{2}$ Université Paris-Est, Cosys, Grettia, Ifsttar, France.}%
\thanks{$^{3}$ LVMT, Enpc, France. } %
\thanks{* Corresponding author. {\texttt{florian.schanzenbacher@ratp.fr}}} }

\newtheorem{theorem}{Theorem}
\newtheorem{corollary}{Corollary}

\graphicspath{{./}}

\begin{document}

\maketitle
\thispagestyle{empty}
\pagestyle{empty}

\begin{abstract}
This paper presents a mathematical model for the
train dynamics in a mass-transit metro line system with one symmetrically
operated junction. We distinguish three parts: a central
part and two branches. The tracks are spatially discretized into
segments (or blocks) and the train dynamics are described by a
discrete event system where the variables are the $k^{th}$ departure
times from each segment. The train dynamics are based on two
main constraints: a travel time constraint modeling theoretic
run and dwell times, and a safe separation constraint modeling
the signaling system in case where the traffic gets very dense.
The Max-plus algebra model allows
to analytically derive the asymptotic average train frequency as a
function of many parameters, including train travel times,
minimum safety intervals, the total number of trains
on the line and the number of trains on each branch. This
derivation permits to understand the physics of traffic. In a further step,
the results will be used for traffic control.

\end{abstract}


\section{Introduction and Literature review}
\label{introduction}

Mass-transit metro lines are generally driven on their capacity limit to satisfy the high demand.
Among the different network topologies for metro lines, those with a junction are
much more sensible to perturbations than simple ring or linear lines.
RATP, the operator of the French capital's metro
system, runs several metro and rapid transit lines with convergences.

The here presented model for a metro line with a symmetrically operated junction describes the
dynamics of the system with static run and dwell times which respect lower
bounds. The model permits a comprehension of the physics of traffic for this case, by
an analytic derivation of the train frequency in function of the parameters of the line
(number of moving trains, safe separation times, etc). 

We base this work on the traffic models presented in~\cite{FNHL16,FNHL16a,FNHL16b}, where the 
physics of traffic in a metro line system without junction are entirely described.
It is a discrete event modeling approach, where we use train departure times as the main
modeling variables. Two cases have been studied in~\cite{FNHL16,FNHL16a,FNHL16b}.

The first case assumes that the train dwell times on platforms respect given lower bounds. 
It has been shown that in this case, the train dynamics can be written linearly in the Max-plus
algebra. This formulation permitted to show that the traffic dynamic admits a unique
asymptotic regime. Moreover, the asymptotic average train time-headway is derived analytically
in function of the number of trains moving on the line and other parameters
like train speed and safe separation time. 




The second case proposes a real-time
control of the train dwell times depending on passenger demand which will be subject to further works on the here presented traffic model.

The Max-plus theory being used here, has further been subject to recent research by Goverde~\cite{Gov07} to analyze railway timetable stability.
Black-box optimization algorithms for real-time railway rescheduling have been treated in research for a long time.
Cacchiani et al.~\cite{Cac14} give a state of the art review of these recovery models and algorithms.
Schanzenbächer et al.~\cite{SCF16} have applied such an optimal holding control for dwell optimization to a mass transit railway line in the Paris area.
Li et al.~\cite{LI17} present an optimal control approach for train regulation and passenger flow control on high-frequency metro lines without a junction.
After a review on Max-plus algebra, we introduce the model of our plant and derive the traffic phases of the system.
\section{Review on linear Max-plus algebra systems}

Max-plus algebra~\cite{BCOQ92} is the idempotent commutative semi-ring $(\mathbb R \cup \{-\infty\}, \oplus, \otimes)$,
where the operations $\oplus$ and $\otimes$ are defined by:
$a \oplus b = \max\{a, b\}$ and $a \otimes b = a + b$. The zero element is $(-\infty)$ denoted by $\varepsilon$ and the
unity element is $0$ denoted by $e$. On the set of square matrices, if $A$ and $B$ are two Max-plus matrices of size
$n \times n$, the addition
$\oplus$ and the product $\otimes$ are defined by: $(A \oplus B)_{ij} := A_{ij} \oplus B_{ij} , \forall i, j$, and
$(AB)_{i,j} = (A \otimes B)_{ij} := \bigoplus_k[A_{ik} \otimes B_{kj}]$. The zero and the unity matrices are also denoted
by $\varepsilon$ and $e$ respectively.

Let us now consider the dynamics of a homogeneous $p$-order Max-plus system with a family of Max-plus matrices $A_l$
\begin{equation}
  x(k) = \bigoplus_{l=0}^p A_l \otimes x(k-l). \label{eq-rev1}
\end{equation}
We define $\gamma$ as the backshift operator applied on the sequences on $\mathbb Z$:
$\gamma^l x(k) = x(k-l), \forall l\in\mathbb N$. Then~(\ref{eq-rev1}) can be written
\begin{equation}
  x(k) = \bigoplus_{l=0}^p \gamma^l A_l x(k) = A(\gamma) x(k),  \label{eq-rev2}
\end{equation}
where $A(\gamma)=\bigoplus_{l=0}^p \gamma^l A_l$ is a polynomial matrix in the backshift operator $\gamma$;
see~\cite{BCOQ92,Gov07} for more details.

$\mu \in \mathbb R_{\max} \setminus \{\varepsilon\}$ is said to be a 
\textit{generalized eigenvalue}~\cite{CCGMQ98} of $A(\gamma)$, with associated
\textit{generalized eigenvector} $v\in \mathbb R_{\max}^n \setminus \{\varepsilon\}$, if
\begin{equation}\label{eq-gev}
  A(\mu^{-1}) \otimes v = v,
\end{equation}
where $A(\mu^{-1})$ is the matrix obtained by evaluating the polynomial matrix $A(\gamma)$ at $\mu^{-1}$.

A directed graph denoted $\mathcal G (A(\gamma))$ can be associated to a dynamic system of type~(\ref{eq-rev2}).
For every $l, 0\leq l\leq p$, an arc $(i,j,l)$ is associated to each non-null ($\neq \varepsilon$) element
$(i,j)$ of Max-plus matrix $A_l$.
A \emph{weight} $W(i,j,l)$ and a \emph{duration} $D(i,j,l)$ are associated to each arc $(i,j,l)$ in the graph,
with $W(i,j,l) = (A_l)_{ij} \neq \varepsilon$ and $D(i,j,l) = l$. 
Similarly, a weight, resp. duration of a cycle (a directed cycle) in the graph is the standard sum of the weights,
resp. durations of all the arcs of the cycle. Finally, the \textit{cycle mean} of a cycle $c$ with a weight $W(c)$ and
a duration $D(c)$ is $W(c)/D(c)$.
A polynomial matrix $A(\gamma)$ is said to be irreducible, if $\mathcal G(A(\gamma))$ is strongly connected.

We recall here a result that we will use in the next sections.
\begin{theorem} \cite[Theorem 3.28]{BCOQ92} \cite[Theorem 1]{Gov07} \label{th-mpa}
  Let $A(\gamma) = \oplus_{l=0}^p A_l\gamma^l$
  be an irreducible polynomial matrix with acyclic
  sub-graph $\mathcal G(A_0)$. Then $A(\gamma)$ has a unique generalized eigenvalue $\mu > \varepsilon$ and finite eigenvectors $v > \varepsilon$
  such that $A(\mu^{-1}) \otimes v = v$, and $\mu$ is equal to the maximum cycle mean of $\mathcal G(A(\gamma))$, given as follows:
  $\mu = \max_{c\in\mathcal C} W(c) / D(c)$,
  where $\mathcal C$ is the set of all elementary cycles in $\mathcal G(A(\gamma))$.
  Moreover, the dynamic system $x(k) = A(\gamma) x(k)$ admits an asymptotic average growth vector (also called cycle time vector here) $\chi$ 
  whose components are all equal to $\mu$.
\end{theorem}

\section{Train dynamics in a metro line system with a junction}


We extend the approach developed in~\cite{FNHL16,FNHL16a,FNHL16b} by modeling a metro line with a junction.
Let us consider a metro line with one junction as shown in Figure~\ref{fig_line} below.
As in~\cite{FNHL16,FNHL16a,FNHL16b}, the line is discretized in a number of segments (or sections, or blocks).
We call node here the point separating two consecutive segments
on the line. Segments and nodes are indexed as in Figure~\ref{fig_line}.

\begin{figure}[h]
 \centering
  \includegraphics[width=0.5\textwidth]{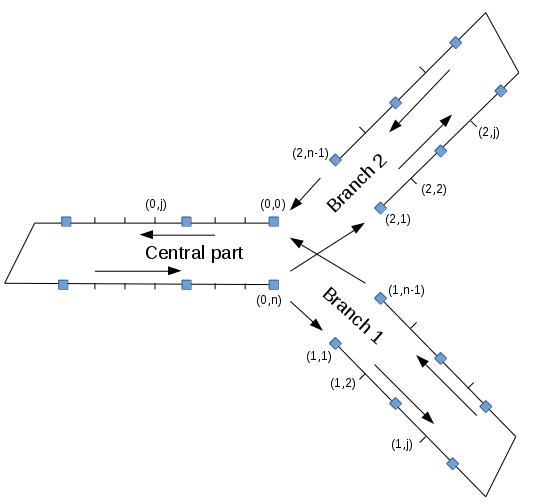} 
\caption {Schema of a metro line with one junction and the corresponding notation.}
\label{fig_line}
\end{figure}

Let us consider the following notations:

\begin{tabular}{ll}
  $u$ & $\in \mathcal U = \{0,1,2\}$ indexes the central part if  \\
      & $u=0$, branch 1 if $u=1$, branch 2 if $u=2$. \\
  $n_u$ & the number of segments on part $u$ of the line. \\
  $m_u$ & the number of trains being on the part $u$ of the  \\
        & line, at time zero. \\
  $b_{(u,j)}$ & $\in \{0,1\}$. It is $O$ (resp. $1$) if there is no \\
              & train (resp. one train) at segment $j$ of part $u$.\\
  $\bar{b}_{(u,j)}$ & $= 1 - b_{(u,j)}$.  \\
  $d^k_{(u,j)}$ & the $k^{\text{th}}$ departure time from node $j$, on part $u$\\
      & of the line. Notice that $k$ does not index trains, \\
      & but count the number of train departures. \\  
  $a^k_{(u,j)}$ & the $k^{\text{th}}$ arrival time to node $j$, on part $u$\\
      & of the line. \\           
  $r{(u,j)}$ & the average running time of trains on segment $j$ \\
      & (between nodes $j-1$ and $j$) of part $u$.\\
  $w^k_{(u,j)}$ & $=d^k_{(u,j)}-a^k_{(u,j)}$ the $k^{\text{th}}$ dwell time on node $j$ \\
                & on part $u$ of the line.\\
  $t^k_{(u,j)}$ & $=r_{(u,j)}+w^k_{(u,j)}$ the $k^{\text{th}}$ travel time from node \\
		& $j-1$ to node $j$ on part $u$ of the line.\\      
  $g^k_{(u,j)}$ & $=a^k_{(u,j)}-d^{k-1}_{(u,j)}$ the $k^{\text{th}}$ safe separation time \\
                & (or close-in time) at node $j$ on part $u$.
\end{tabular}  
\begin{tabular}{ll}      
  $h^k_{(u,j)}$ & $=d^k_{(u,j)}-d^{k-1}_{(u,j)} = g^k_{(u,j)}+w^{k}_{(u,j)}$ the $k^{\text{th}}$\\
      & departure time-headway at node $j$ on part $u$.\\
  $s^k_{(u,j)}$ & $=g^{k+b_{(u,j)}}_{(u,j)}-r_{(u,j)}$.    
\end{tabular}

~

We also use underline notations to note the minimum 
bounds of the corresponding variables respectively.
Then $\underline{r}_{(u,j)}, \underline{t}_{(u,j)}, \underline{w}_{(u,j)}, \underline{g}_{(u,j)}, 
\underline{h}_{(u,j)}$ and $\underline{s}_{(u,j)}$
denote respectively minimum running, travel, dwell, safe separation, headway and 
$s$ times.

The (asymptotic) averages on $j$ and on $k$ of those variables are denoted 
without any subscript or superscript.
Then $r, t, w, g, h$ and $s$ denote the average running, travel, dwell, safe separation,
 headway and $s$ times, respectively.

It is easy to check the following relationships:  
  \begin{align}
    & g_u = r_u + s_u, \label{form1} \\
    & t_u = r_u + w_u, \label{form2} \\
    & h_u = g_u + w_u = r_u + w_u + s_u = t_u + s_u.
    \label{form3}
  \end{align}    
We give below the train dynamics in the metro line system. We distinguish the dynamics on the tracks out
of the junction, with the ones on the divergence and on the merge.

\subsection{Train dynamics out of the junction}

We model the train dynamics here as in~\cite{FNHL16,FNHL16a,FNHL16b}.
Two main constraints are considered to describe the train dynamics out of the junction.
\begin{itemize}
  \item The $k^{\text{th}}$ train departure from node $j$ (of any part of the line) corresponds 
     to the $k^{\text{th}}$ train departure from node $j-1$ in case where there was no train 
     on segment $j$ at time zero; and it corresponds to the $(k-1)^{\text{th}}$ 
     train departure from node $j-1$ in case where there was a train on segment $j$ at time zero.
     Between two consecutive departures, a minimum time of $\underline{t}_{(u,j)}$ is respected.
     We write
     \begin{equation}\label{eq-c1}
        d^k_{(u,j)} \geq d^{k-b_{(u,j)}}_{(u,j-1)} + \underline{t}_{(u,j)}, \; \forall k\geq 0, u\in\mathcal U, j \neq n_u.
     \end{equation}
   
  \item The $k^{\text{th}}$ train departure from node $j$ must be preceded by the $(k-1)^{\text{th}}$ train departure
    from node $j+1$ plus a minimum time $\underline{s}^k_{(u,j)}$ in case where there was no train on segment $j+1$ at time zero; and
    it must be preceded by the $k^{\text{th}}$ train departure from node $j+1$ plus a minimum time $\underline{s}^k_{(u,j)}$
    in case where there was a train on segment $j+1$ at time zero. We write
    \begin{equation}\label{eq-c2}
          d^k_{(u,j)} \geq d^{k-\bar{b}_{(u,j+1)}} _{(u,j+1)} + \underline{s}_{(u,j+1)}, \; \forall k\geq 0, u\in\mathcal U, j \neq n_u.
    \end{equation}
\end{itemize}

We assume here that a train departs from node $j$ out of the junction, as soon as the two
constraints~(\ref{eq-c1}) and~(\ref{eq-c2}) are satisfied. We get


\begin{equation}\label{eq-c3}
  d^k_{(u,j)} = \max \left\{ d^{k-b_{(u,j)}}_{(u,j-1)} + \underline{t}_{(u,j)}, d^{k-\bar{b}_{(u,j+1)}} _{(u,j+1)} + \underline{s}_{(u,j+1)}  \right\}.
\end{equation}
This assumption holds for all couples of constraints we will propose below.
This will permit us to write the whole train dynamics as a homogeneous Max-plus system.

\subsection{Train dynamics on the divergence}

We assume here that trains leaving the central part of the line and 
going to the branches respect the following rule.
Odd departures go to branch 1 while even departures go to branch~2.
We then have the following constraints.

The $k^{\text{th}}$ departures from the central part:
\begin{equation}\label{eq-d1}
  d^k_{(0,n)} \geq d^{k-b_{(0,n)}}_{(0,n-1)} + \underline{t}_{(0,n)}, \; \forall k\geq 0, \\
\end{equation}

\begin{equation}\label{eq-d2}
  d^k_{(0,n)} \geq \begin{cases}
                      d^{(k+1)/2-\bar{b}_{(1,1)}}_{(1,1)} + \underline{s}_{(1,1)} & \text{if } k \text{ is odd} \\ ~~ \\
                      d^{k/2-\bar{b}_{(2,1)}}_{(2,1)} + \underline{s}_{(2,1)} & \text{if } k \text{ is even} \\
                   \end{cases}
\end{equation}

The $k^{\text{th}}$ departures from the entry of branch 1:
\begin{eqnarray}
   d^k_{(1,1)} \geq d^{(2k-1)-b_{(1,1)}}_{(0,n)} + \underline{t}_{(1,1)}, \; \forall k\geq 0, \label{eq-d3}\\
   d^k_{(1,1)} \geq d^{k-\bar{b}_{(1,2)}} _{(1,2)} + \underline{s}_{(1,2)}, \; \forall k\geq 0. \label{eq-d4}
\end{eqnarray}

The $k^{\text{th}}$ departures from the entry of branch 2:
\begin{eqnarray}
   d^k_{(2,1)} \geq d^{2k-b_{(2,1)}}_{(0,n)} + \underline{t}_{(2,1)}, \; \forall k\geq 0, \label{eq-d5}\\
   d^k_{(2,1)} \geq d^{k-\bar{b}_{(2,2)}}_{(2,2)} + \underline{s}_{(2,2)}, \; \forall k\geq 0. \label{eq-d6}
\end{eqnarray}

\subsection{Train dynamics on the merge}

We assume here that trains entering to the central part of the line from 
the two branches respect the following rule.
Odd departures at node $(0,0)$ towards the central part correspond to trains coming from
branch~1 while even ones correspond to trains coming from branch~2.

The $k^{\text{th}}$ departures from the central part:
\begin{equation}\label{eq-m1}
  d^k_{(0,0)} \geq \begin{cases}
                      d^{(k+1)/2-b_{(1,n)}}_{(1,n-1)} + \underline{t}_{(1,n)} & \text{if } k \text{ is odd} \\ ~~ \\
                      d^{k/2-b_{(2,n)}}_{(2,n-1)} + \underline{t}_{(2,n)} & \text{if } k \text{ is even} \\
                   \end{cases}
\end{equation}
\begin{equation}\label{eq-m2}
  d^k_{(0,0)} \geq d^{k-\bar{b}_{(0,1)}}_{(0,1)} + \underline{s}_{(0,1)}, \; \forall k\geq 0, \\
\end{equation}

The $k^{\text{th}}$ departures from the entry of branch 1:
\begin{eqnarray}
   d^k_{(1,n-1)} \geq d^{k-b_{(1,n-1)}}_{(1,n-2)} + \underline{t}_{(1,n-1)}, \; \forall k\geq 0, \label{eq-m3}\\
   d^k_{(1,n-1)} \geq d^{(2k-1)-\bar{b}_{(1,n)}}_{(0,0)} + \underline{s}_{(1,n)}, \; \forall k\geq 0. \label{eq-m4}
\end{eqnarray}

The $k^{\text{th}}$ departures from the entry of branch 2:
\begin{eqnarray}
   d^k_{(2,n-1)} \geq d^{k-b_{(2,n-1)}}_{(2,n-2)} + \underline{t}_{(2,n-1)}, \; \forall k\geq 0, \label{eq-m5} \\
   d^k_{(2,n-1)} \geq d^{2k-\bar{b}_{(2,n)}}_{(0,0)} + \underline{s}_{(2,n)}, \; \forall k\geq 0. \label{eq-m6}
\end{eqnarray}

\section{Max-plus algebra modeling}

In order to avoid multiplicative backshifts between the departures
on the central part and the ones on the branches, we introduce a
change of variables in this section. 
Let us look at the dynamic given by~(\ref{eq-d3}).
$d^k_{(1,1)}$ is given as a function of $d^{2k-1-b_{(1,1)}}_{(0,n)}$.
We see clearly that the two sequences do not have the same growth speed.
Indeed, the growth rate of $d_{(0,n)}$ is about double the one of $d_{(1,1)}$.
This is due to the fact that every second train moving on the central part of the line
goes to branch~1. 

In order to have all the sequences of the train dynamics growing with the same speed,
and then be able to write the dynamics as a homogeneous Max-plus system, we
introduce the following change of variables:
\begin{align}
  & \delta^k_{(0,j)} = d^k_{(0,j)}, \forall k\geq 0, \forall j  \\
  & \delta^{2k}_{(1,j)} = d^k_{(1,j)}, \forall k\geq 0, \forall j  \\
  & \delta^{2k}_{(2,j)} = d^k_{(2,j)}, \forall k\geq 0, \forall j.
\end{align}
In the following, we rewrite all the train dynamics with the change 
of variables defined above.

\subsection{New train dynamics out of the junction}

The train dynamics out of the junction are written as follows.

On the central part, it is sufficient to replace $d$ with $\delta$:
\begin{eqnarray}
   \delta^k_{(0,j)} \geq \delta^{k-b_{(0,j)}}_{(0,j-1)} + \underline{t}_{(0,j)}, \; \forall k\geq 0, j \neq n_u, \label{eq-nc1}\\
   \delta^k_{(0,j)} \geq \delta^{k-\bar{b}_{(0,j+1)}} _{(0,j+1)} + \underline{s}_{(0,j+1)}, \; \forall k\geq 0, j \neq n_u. \label{eq-nc2}
\end{eqnarray}

For the dynamics on the branches, we get
\begin{equation}
   \delta^k_{(u,j)} \geq \delta^{k-2b_{(u,j)}}_{(u,j-1)} + \underline{t}_{(u,j)}, \; \forall k\geq 0, u\in\{1,2\}, j \neq n_u, \label{eq-nc3}
\end{equation}
\begin{equation}
   \delta^k_{(u,j)} \geq \delta^{k-2\bar{b}_{(u,j+1)}} _{(u,j+1)} + \underline{s}_{(u,j+1)}, \; \forall k\geq 0, u\in\{1,2\}, j \neq n_u. \label{eq-nc4}
\end{equation}

\subsection{New train dynamics on the divergence}

The dynamics on the divergence are rewritten as follows.

The $k^{\text{th}}$ departures from the central part:
\begin{equation}\label{eq-nd1}
  \delta^{k}_{(0,n)} \geq \delta^{k-b_{(0,n)}}_{(0,n-1)} + \underline{t}_{(0,n)}, \; \forall k\geq 0, \\
\end{equation}

\begin{equation}\label{eq-nd2}
  \delta^{k}_{(0,n)} \geq \begin{cases}
                      \delta^{(k+1)-2\bar{b}_{(1,1)}}_{(1,1)} + \underline{s}_{(1,1)} & \text{for } k \text{ is odd} \\ ~~ \\
                      \delta^{k-2\bar{b}_{(2,1)}}_{(2,1)} + \underline{s}_{(2,1)} & \text{for } k \text{ is even} \\
                   \end{cases}
\end{equation}

The $k^{\text{th}}$ departures from the entry of branch 1:
\begin{eqnarray}
   \delta^{k}_{(1,1)} \geq \delta^{(k-1)-2b_{(1,1)}}_{(0,n)} + \underline{t}_{(1,1)}, \; \forall k\geq 0, \label{eq-nd2}\\
   \delta^{k}_{(1,1)} \geq \delta^{k-2\bar{b}_{(1,2)}} _{(1,2)} + \underline{s}_{(1,2)}, \; \forall k\geq 0. \label{eq-n3}
\end{eqnarray}

The $k^{\text{th}}$ departures from the entry of branch 2:
\begin{eqnarray}
   \delta^k_{(2,1)} \geq \delta^{k-2b_{(2,1)}}_{(0,n)} + \underline{t}_{(2,1)}, \; \forall k\geq 0, \label{eq-nd4} \\
   \delta^k_{(2,1)} \geq \delta^{k-2\bar{b}_{(2,2)}}_{(2,2)} + \underline{s}_{(2,2)}, \; \forall k\geq 0. \label{eq-nd5}
\end{eqnarray}

\subsection{New train dynamics on the merge}

The dynamics on the merge are rewritten as follows.

The $k^{\text{th}}$ departures from the central part:
\begin{equation} \label{eq-nm1}
  \delta^k_{(0,0)} \geq \begin{cases}
                      \delta^{(k+1)-2b_{(1,n)}}_{(1,n-1)} + \underline{t}_{(1,n)} & \text{for } k \text{ is odd} \\ ~~ \\
                      \delta^{k-2b_{(2,n)}}_{(2,n-1)} + \underline{t}_{(2,n)} & \text{for } k \text{ is even} \\
                   \end{cases}
\end{equation}
\begin{equation}\label{eq-nm2}
  \delta^k_{(0,0)} \geq \delta^{k-2\bar{b}_{(0,1)}}_{(0,1)} + \underline{s}_{(0,1)}, \; \forall k\geq 0, \\
\end{equation}

The $k^{\text{th}}$ departures from the entry of branch 1:
\begin{eqnarray}
   \delta^k_{(1,n-1)} \geq \delta^{k-2b_{(1,n-1)}}_{(1,n-2)} + \underline{t}_{(1,n-1)}, \; \forall k\geq 0, \label{eq-nm3} \\
   \delta^k_{(1,n-1)} \geq \delta^{(k-1)-2\bar{b}_{(1,n)}}_{(0,0)} + \underline{s}_{(1,n)}, \; \forall k\geq 0. \label{eq-nm4}
\end{eqnarray}

The $k^{\text{th}}$ departures from the entry of branch 2:
\begin{eqnarray}
   \delta^k_{(2,n-1)} \geq \delta^{k-2b_{(2,n-1)}}_{(2,n-2)} + \underline{t}_{(2,n-1)}, \; \forall k\geq 0, \label{eq-nm5} \\
   \delta^k_{(2,n-1)} \geq \delta^{k-2\bar{b}_{(2,n)}}_{(0,0)} + \underline{s}_{(2,n)}, \; \forall k\geq 0. \label{eq-nm6}
\end{eqnarray}

\subsection{Train dynamics in Max-plus algebra}

Let us now show how all the train dynamics given above are written in Max-plus algebra.
First, as already mentioned above, we assume for every couple of constraints 
written on the departure from a given node, that the departure in question is realized as
soon as the two associated constraints are satisfied.
For example, with this assumption, constraints~(\ref{eq-nc1}) and~(\ref{eq-nc2}) give
\begin{equation}\label{eq-mp1}
   \delta^k_{(0,j)} = \max \left\{ \delta^{k-b_{(0,j)}}_{(0,j-1)} + \underline{t}_{(0,j)},  \delta^{k-\bar{b}_{(0,j+1)}} _{(0,j+1)} + \underline{s}_{(0,j+1)} \right\}.
\end{equation}
which is written in Max-plus algebra as follows:
\begin{equation}\label{eq-mp2}
   \delta_{(0,j)} = \gamma^{b_{(0,j)}} \underline{t}_{(0,j)} \delta_{(0,j-1)} \oplus  \gamma^{\bar{b}_{(0,j+1)}} \underline{s}_{(0,j+1)} \delta_{(0,j+1)}.
\end{equation}

We can easily check that all the couples of constraints of the whole dynamics can now (with the change of variables) be written in Max-plus algebra,
as done above in~(\ref{eq-mp2}).
However, because of the junction, where \textit{every} $2^{nd}$ \textit{train} goes in the alternative direction,
we will get two different homogeneous Max-plus systems that are applied 
alternatively for odd and even $k^{th}$ \textit{departures}. 
If we denote $\delta^k = {}^t (\delta_0^k,\delta_1^k,\delta_2^k)$
the column vector which concatenates the three vectors $\delta_0^k, \delta_1^k$
and $\delta_2^k$ (with $\delta_0^k$ the colum vector with components $\delta^k_{(0,j)}$), then the whole train dynamics can be written as follows:
\begin{equation}\label{eq-mp3}
  \delta^k =   \begin{cases}
		  A^{(1)}(\gamma) \otimes \delta^k & \text{ if } k \text{ is odd},\\
		  A^{(2)}(\gamma) \otimes \delta^k & \text{ if } k \text{ is even},
	       \end{cases}
\end{equation}
where 
\begin{equation}
   A^{(1)}(\gamma) = \begin{pmatrix}
                A^{(1)}_{00}(\gamma) & A^{(1)}_{01}(\gamma) & \varepsilon \\
                A^{(1)}_{10}(\gamma) & A^{(1)}_{11}(\gamma) & \varepsilon \\
                \varepsilon & \varepsilon & A^{(1)}_{22}(\gamma)
             \end{pmatrix},
\end{equation}
\begin{equation}
   A^{(2)}(\gamma) = \begin{pmatrix}
                A^{(2)}_{00}(\gamma) & \varepsilon & A^{(2)}_{02}(\gamma) \\
                \varepsilon & A^{(2)}_{11}(\gamma) & \varepsilon \\
                A^{(2)}_{20}(\gamma) & \varepsilon & A^{(2)}_{22}(\gamma)
             \end{pmatrix}.
\end{equation}

The diagonal blocks of the matrices above are given as follows 
($\forall u\in \{0,1,2\}$ and $p\in\{1,2\}$ and for $n_u=10$ as an example):

\small
$A^{(p)}_{uu}(\gamma) = $
$$		 \begin{pmatrix}
                    \varepsilon & \gamma^{\bar{b}_{(u,2)}} \underline{s}_{(u,2)} & \ldots & \varepsilon \\
                    \gamma^{b_{(u,2)}} \underline{t}_{(u,2)} & \varepsilon & \ddots & \varepsilon  \\
                    \varepsilon & \ddots & \ddots & \gamma^{\bar{b}_{(u,9)}} \underline{s}_{(u,9)} \\
                    \varepsilon & \varepsilon & \gamma^{b_{(u,8)}} \underline{t}_{(u,8)} & \varepsilon
                 \end{pmatrix}.$$
\normalsize

To have an idea of the other blocks we give here $A^{(1)}_{01}(\gamma)$:

\small
$$ A^{(1)}_{01}(\gamma) = \begin{pmatrix}
                    \varepsilon & \varepsilon & \cdots &\gamma^{b_1} \underline{t}_1 \\
                    \varepsilon & \varepsilon & \cdots & \varepsilon  \\
                    \vdots &  & \ddots & \vdots \\
                    \gamma^{\bar{b}_{1,j=6}} \underline{s}_{1,j=5} & \varepsilon & \cdots & \varepsilon
                 \end{pmatrix}.$$
\normalsize

We keep in mind that by the changing of variables done above, the number of $k^{th}$ \textit{departures} on the branches
has been doubled. Most importantly, this means that we have one average asymptotic growth rate for the matrices
$A^{(1)}, A^{(2)}$.


To correctly represent the junction, we consider the composition 
of the train dynamics with itself, which gives us the dynamics on two steps.
We get matrix \textit{B}, whose average asymptotic growth rate is equal to the average time-headway
between two consecutive $k^{th}$ \textit{departures} 
(e.g. the time-headway between two trains going in different directions), and therefore represents
the average time-headway on the central part.
\begin{equation}\label{eq-mp4}
  \delta^k =  B(\gamma) \otimes \delta^k,
\end{equation}
where $B(\gamma) = A^{(2)}(\gamma) \otimes A^{(1)}(\gamma)$.

\section{Analytical derivation of traffic phases}

We will now show how the Max-plus model allows to derive the average train time-headway.
We present the results of an application to a metro line with a junction in Paris, France, and compare analytical results
to simulation and to the actual timetable.
Let us notice that if the growth rate $h$ of system~(\ref{eq-mp4}) exists, it represents the time-headway on the central part, and
since the number of \textit{k steps} on the branches has been doubled because of the changing of variables,
the time-headway on the branches is $2h$.
The growth rate is given by the unique generalized eigenvalue of the homogeneous Max-plus system, which can be calculated
from its associated graph (Theorem~\ref{thm-1}).



We show that the asymptotic average 
train frequency of a metro line with a junction, depends on the total number of trains and on the difference between the number
of trains on the branches. Both parameters are invariable in time (in two steps of the train dynamics),
since the rule \textit{every} $2^{nd}$ \textit{train} is applied on the divergence and on the merge.
We consider the following notations:

\begin{tabular}{ll}
  $m$ & $= m_0+m_1+m_2$ the total number of trains \\
      & on the line.\\
  $\Delta m$ & $= m_2 - m_1$ the difference in the number \\
     & of trains between branches 2 and 1.\\
  $\bar{m}_u$ & $= n_u - m_u, \forall u\in\{0,1,2\}$.\\      
  $\bar{m}$ & $= \bar{m}_0 + \bar{m}_1 + \bar{m}_2$.\\
  $\Delta \bar{m}$ & $= \bar{m}_2 - \bar{m}_1$.\\
  $\underline{T}_u$ & $= \sum_j \underline{t}_{(u,j)}, \forall u\in\{0,1,2\}$. \\
  $\underline{S}_u$ & $= \sum_j \underline{s}_{(u,j)}, \forall u\in\{0,1,2\}$.
\end{tabular}

~~

\begin{theorem}\label{thm-1}
  The dynamic system~(\ref{eq-mp4}) admits a unique
  asymptotic stationary regime, with a common average growth
  rate $h_0$ for all the variables, which represents the average train
  time-headway $h_0$ on the central part and $h_1/2 = h_2/2$ on the branches.
  Moreover we have
  $$ h_0 = h_1/2 = h_2/2 = \max \{ h_{fw}, h_{\min} ,h_{bw}, h_{br}\}, $$
  with\footnote{fw: forward, bw: backward, min: minimum, br: branches.}
  $$ h_{fw} = \max\left\{ \frac{\underline{T}_0 + \underline{T}_1}{m - \Delta m},
                          \frac{\underline{T}_0 + \underline{T}_2}{m + \Delta m} \right\}, $$
	$$ h_{\min} = \max \begin{cases}
		  \max_{u,j} (t_{(u,j)} + s_{(u,j)}) & \forall u \in \{0\},\\
		  \max_{u,j} (t_{(u,j)} + s_{(u,j)})/2 & \forall u \in \{1,2\},
	\end{cases}$$
  $$ h_{bw} = \max\left\{ \frac{\underline{S}_0 + \underline{S}_1}{\bar{m} - \Delta \bar{m}},
                          \frac{\underline{S}_0 + \underline{S}_2}{\bar{m} + \Delta \bar{m}} \right\}, $$
  $$ h_{br} =  \max\left\{\frac{\underline{T}_1 + \underline{S}_2}{2(n_2 - \Delta m)},
                  \frac{\underline{S}_1 + \underline{T}_2}{2(n_1 + \Delta m)}\right\}.$$
\end{theorem}

~~

\begin{figure}[h]
\centering
\tikzset{
    side by side/.style 2 args={
        line width=1.0pt,
        #1,
        postaction={
            clip,postaction={draw,#2}
        }
    }
}
\tikzset{
    side by side 2/.style 2 args={
        line width=1.5pt,
        #1,
        postaction={
            clip,postaction={draw,#2}
        }
    }
}
\definecolor{light-gray}{gray}{0.75}
\begin{tikzpicture}[->,auto, node distance=1.9 cm, state_0/.style={circle,fill=black,draw,text=white},state_1/.style={circle,fill=light-gray,draw,text=white}
,state_2/.style={circle,fill=gray,draw,text=white}]
\node [state_0] (1) {0,0};
\node [state_0](2) [above left of=1] {0,1};
\node [state_0](4) [below of=1] {0,3};
\node [state_0](3) [below left of=4] {0,2};
\node [state_1](5) [below right of=4] {1,1};
\node [state_1](6) [below right of=5] {1,2};
\node [state_1](7) [above right of=6] {1,3};
\node [state_1](8) [above right of=5] {1,4};
\node [state_2](9) [above of=8] {2,1};
\node [state_2](10) [above right of=9] {2,2};
\node [state_2](11) [above left of=10] {2,3};
\node [state_2](12) [below left of=11] {2,4};
\path[]
(1) edge [line width=1 pt,red] node [] {} (3)
(2) edge [line width=1 pt,red] node [] {} (4)
(3) edge [bend left,line width=1 pt,blue] node [] {} (1)
(4) edge [bend left,line width=1 pt,blue] node [] {} (2)
(5) edge [bend left,line width=1 pt,red] node [] {} (7)
(6) edge [bend right,line width=1 pt,red] node [] {} (8)
(7) edge [bend left,line width=1 pt,blue] node [] {} (5)
(8) edge [bend right,line width=1 pt,blue] node [] {} (6)
(9) edge [bend right,line width=1 pt,red] node [] {} (11)
(10) edge [bend left,line width=1 pt,red] node [] {} (12)
(11) edge [bend right,line width=1 pt,blue] node [] {} (9)
(12) edge [bend left,line width=1 pt,blue] node [] {} (10)
(4) edge [bend right,line width=1 pt,red] node [] {} (6)
(8) edge [bend right,line width=1 pt,red] node [] {} (2)
(3) edge [bend right,line width=1 pt,red] node [] {} (9)
(11) edge [bend right,line width=1 pt,red] node [] {} (1)
(2) edge [bend left,line width=1 pt,blue] node [] {} (12)
(10) edge [bend left,line width=1 pt,blue] node [] {} (4)
(1) edge [bend left,line width=1 pt,blue] node [] {} (7)
(5) edge [bend left,line width=1 pt,blue] node [] {} (3)
(1) edge [loop above,side by side={red}{blue}] node [] {} (1)
(2) edge [loop above,side by side={red}{blue}] node [] {} (2)
(3) edge [loop below,side by side={red}{blue}] node [] {} (3)
(4) edge [loop below,side by side={red}{blue}] node [] {} (4)
(5) edge [loop below,side by side={red}{blue}] node [] {} (5)
(6) edge [loop right,side by side={red}{blue}] node [] {} (6)
(7) edge [loop above,side by side={red}{blue}] node [] {} (7)
(8) edge [loop right,side by side={red}{blue}] node [] {} (8)
(9) edge [loop right,side by side={red}{blue}] node [] {} (9)
(10) edge [loop below,side by side={red}{blue}] node [] {} (10)
(11) edge [loop right,side by side={red}{blue}] node [] {} (11)
(12) edge [loop above,side by side={red}{blue}] node [] {} (12)
(8) edge [bend left,line width=1.5 pt,side by side ={red}{blue}] node [] {} (12)
(12) edge [right,line width=0.5 pt,gray] node [] {} (8)
(5) edge [bend left,line width=1.5 pt,side by side ={red}{blue}] node [] {} (9)
(9) edge [right,line width=0.5 pt,gray] node [] {} (5);
\begin{customlegend}[legend cell align=left, 
legend entries={ 
$W=t_{(u,j)}$: travel time,
$W=s_{(u,j)}$: s time,
\textit{(u,j) = (0,j)}: node on central part, 
\textit{(u,j) = (1,j)}: node on branch 1,
\textit{(u,j) = (2,j)}: node on branch 2
},
legend style={at={(1.9,-4.4)},font=\footnotesize}] 
    \addlegendimage{-stealth,red,opacity=1}
    \addlegendimage{-stealth,blue,opacity=1}
    \addlegendimage{mark=ball,ball color = black,draw=white}
    \addlegendimage{mark=ball,ball color=light-gray,draw=white}
    \addlegendimage{mark=ball,ball color=gray,draw=white}
\end{customlegend}
\end{tikzpicture}
\caption{$\mathcal G(B(\gamma))$ for $n_0=3$, $n_1=n_2=5$.}
\label{fig-ag}
\end{figure}
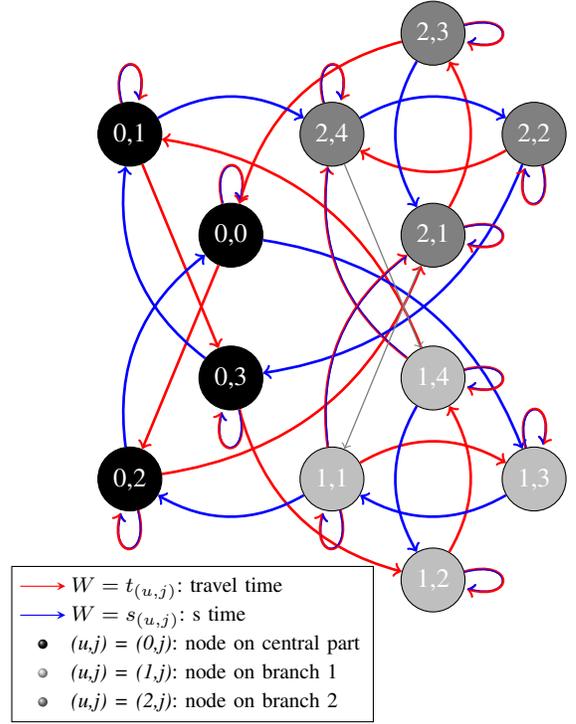

\proof 
It consists in applying Theorem~\ref{th-mpa}, which gives $h_0$ as 
the maximum cycle mean of $\mathcal G(B(\gamma))$.
We give here (Figure~\ref{fig-ag}) the result for $n_0=3$, $n_1=n_2=5$ (same
type of cycles for any values of $n_u$).


\begin{itemize}
  \item The two red cycles $(0,1)-(0,3)-(1,2)-(1,4)$ and $(0,0)-(0,2)-(2,1)-(2,3)$ in the travel direction.
  \item[$\Rightarrow$] $h_{fw}$ = maximum of the cycle means of these two cycles.
  \item The two blue cycles $(0,1)-(2,4)-(2,2)-(0,3)$ and $(0,0)-(1,3)-(1,1)-(0,2)$ against the travel direction.
  \item[$\Rightarrow$] $h_{bw}$ = maximum of the cycle means of these two cycles.
  \item The red/blue loops on all the nodes of the graph.
  \item[$\Rightarrow$] $h_{\min}$ = maximum of the cycle means of the loops.
  \item The cycles with two arcs $(0,0)-(0,2)$, $(0,1)-(0,3)$, $(1,1)-(1,3)$, 
    $(1,2)-(1,4)$, $(2,1)-(2,3)$, $(2,2)-(2,4)$.
  \item[$\Rightarrow$] Their cycle means are dominated by those of the loops, since their mean
    is the average of two neighbored loops.
\end{itemize}
\begin{itemize}
  \item The two cycles $(0,0)-(1,3)-(1,1)-(2,1)-(2,3)$ and $(0,3)-(1,2)-(1,4)-(2,4)-(2,2)$ passing
    by the two branches, one in the travel direction, the other against the travel direction.
  \item[$\Rightarrow$] $h_{br}$ = maximum of the cycle means of these two cycles.
\end{itemize}
\endproof

\noindent
\begin{figure*}[t]
  \centering
  \begin{tabular}{cc}
     \centering
     \includegraphics[scale=0.23]{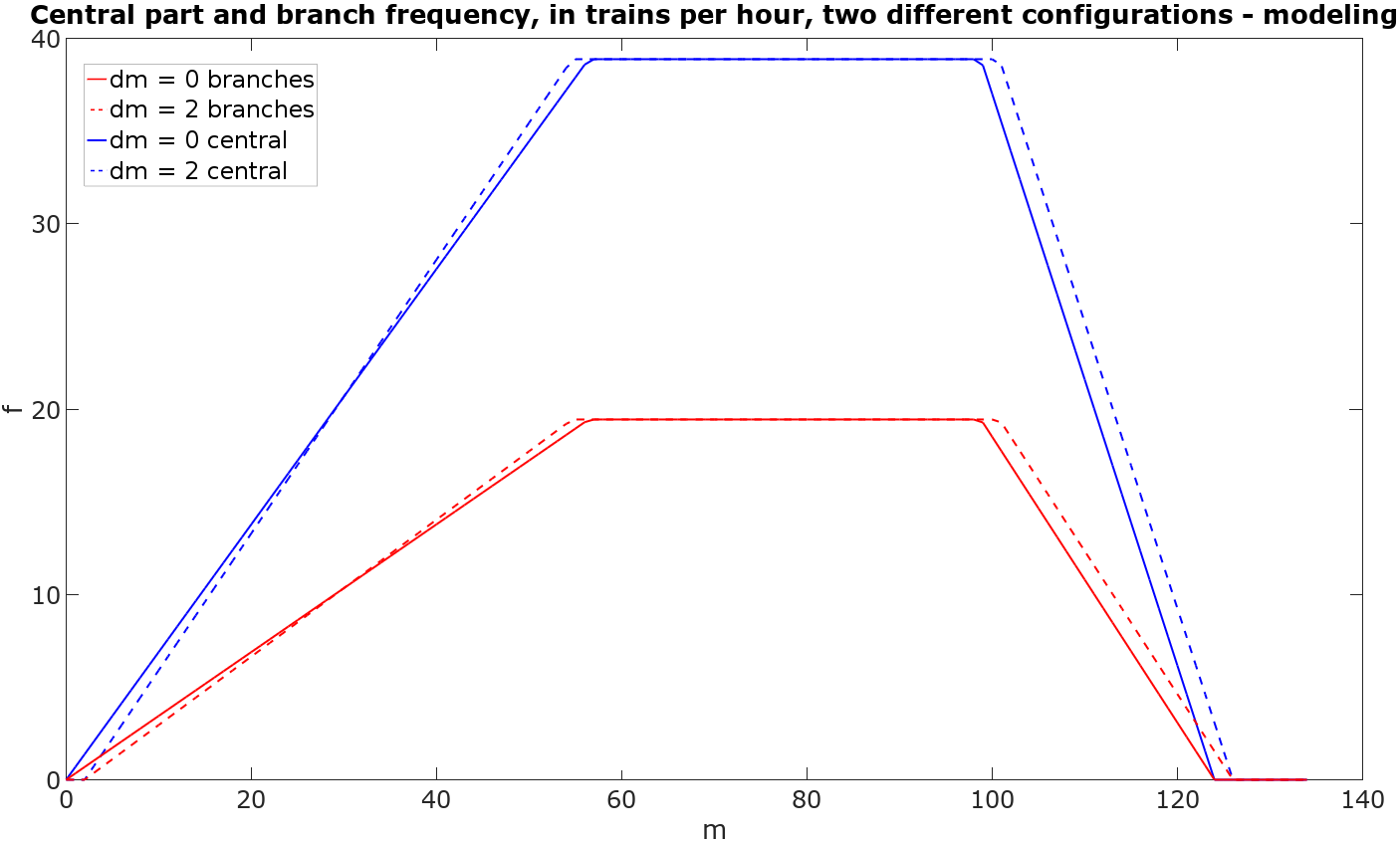} & \includegraphics[scale=0.23]{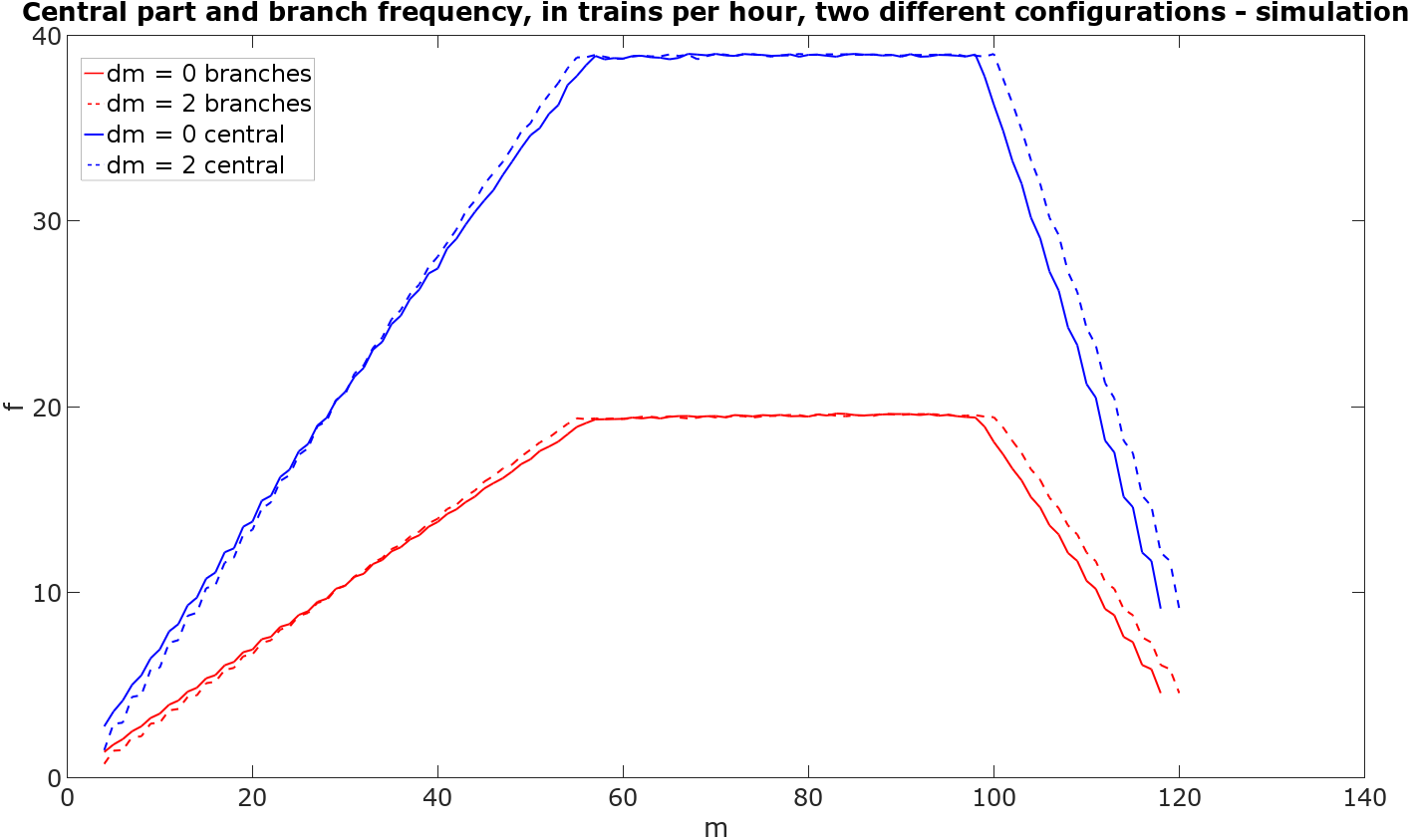}
  \end{tabular}
  \caption{The asymptotic average train frequency $f$ (blue: central part, red: branches), displayed as a function of the number $m$ of moving trains, for
     the two cases of $\Delta m = 0$ (solid line) and $\Delta m = 2$ (dashed line). On the left side: analytically derived formula (Corollary~\ref{cor1}).
     On the right side: simulation.}
  \label{fig-mod/sim}
\end{figure*}

\begin{corollary}\label{cor1}
  The average train frequency $f_0$ on the central part and $f_1=f_2$ on the branches
  are given as follows:
  $$f_0 = 2 f_1 = 2 f_2 = \max \left\{ 0, \min \left\{ \frac{1}{h_{fw}}, \frac{1}{h_{\min}} ,\frac{1}{h_{bw}}, \frac{1}{h_{br}} \right\} \right\}.$$
\end{corollary}

~~

\proof
Directly from Theorem~\ref{thm-1}, with $0\leq f=1/h$.
\endproof




Theorem~\ref{thm-1} shows that in a metro line system with a junction 
and two symmetrically operated branches, the part with the longest
time-headway imposes its frequency to the rest of the system (with the
frequency on the branches being half the one on the central part).

\begin{figure}[h]
  \centering
  \includegraphics[scale=0.24]{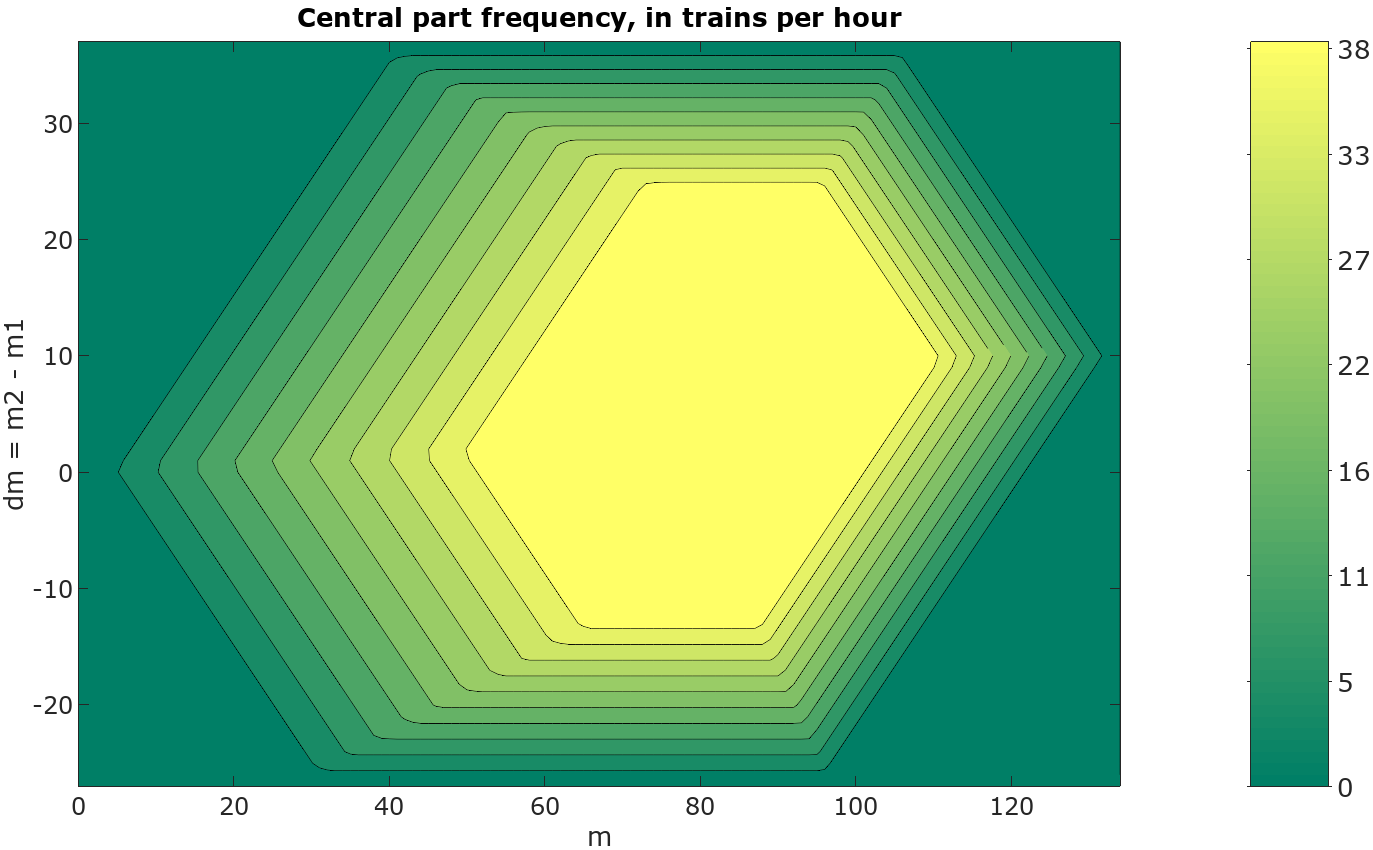}
  \caption{The asymptotic average train frequency $f_0$, as a function of the number $m$ of trains on the line 
          and of the difference $\Delta m$ of the number of trains on the two branches.
          We recognize the eight traffic phases of Theorem~\ref{thm-1} and Corollary~\ref{cor1}.}
  \label{fig-2D}
\end{figure}

We depict in Figure~\ref{fig-2D} the analytically derived traffic phases of the train dynamics.
These frequencies are piecewise linear (Theorem~\ref{thm-1} and Corollary~\ref{cor1}).
RATP, the metro operator of Paris, France, has provided the real values of the minimum running, dwell and safe separation times of a metro line with a junction.
Eight traffic phases can be distinguished. The frequencies of the traffic phases in Figure~\ref{fig-2D}, represent the central part of the line.
A detailed explanation of the phases will be given in a further paper.

We can see, that for every $m$, it exists a $\Delta m$ which maximizes the frequency (Theorem~\ref{thm-2}):

\begin{theorem}\label{thm-2}
$\forall m, \exists \Delta m, f(m, \Delta m) \geq f (m, \Delta m'), \forall \Delta m'$.
\end{theorem}

Theorem~\ref{thm-2} is an important result and will be used in our further research for traffic control.

Figure~\ref{fig-mod/sim} illustrates the traffic phases, derived by Theorem~\ref{thm-1} and Corollary~\ref{cor1}, for two different values of $\Delta m$ on the studied metro line with a junction in Paris, France.
On the left side, the analytically derived results are given. On the right side, we show the results from numerical simulations, for comparison.
Notice that the analytical derivation and the simulation are coherent.
Furthermore, for $m = 52$ and $\Delta m = 2$, the time-headway (and the frequency) of our model represents precisely the timetable of the line.

To illustrate the impact of the parameter $\Delta m$ on the average asymptotic frequency, we give another configuration, $\Delta m = 0$.
Let us notice, that on this line, for $m = 52$, $\Delta m = 2$ maximizes the average asymptotic frequency accordingly to Theorem~\ref{thm-2} (proof not given here).

\section{Conclusion and future work}

These first results of our Max-plus approach to model the dynamical behavior in a metro line system with a junction
are encouraging. We will further develop the model towards dynamic dwell times in order to take into account the
passenger demand on the platforms and in the trains, as well as dynamic running times to recover perturbations and to
stabilize the system.
Finally, our future work will focus on a real-time version of the model, where the system is optimized under dynamic passenger demand to guarantee stability.

%




\end{document}